\documentclass{article}
\usepackage[T1]{fontenc}
\usepackage{amsthm,amsmath,amssymb}
\usepackage{bm}
\usepackage{mathrsfs}
\usepackage[all]{xy}
\usepackage{enumerate}

\newtheorem{theorem}{Theorem}[section]
\newtheorem{proposition}[theorem]{Proposition}
\newtheorem{lemma}[theorem]{Lemma}
\theoremstyle{definition}
\newtheorem{definition}[theorem]{Definition}

\begin{document}

  \title{Natural density and probability, constructively}
\author{Samuele Maschio}
\date{}
\maketitle
\begin{abstract}
  We give here a constructive account of the frequentist approach to probability, by means of natural density. Then we discuss some probabilistic variants of the Limited Principle of Omniscience. \end{abstract}

\section{Introduction}
\emph{Natural density} (\cite{niven1951}) provides a notion of size for subsets of natural numbers. Classically the density of $A\subseteq \mathbb{N}^{+}$ is defined as $$\delta(A)=\lim_{n\rightarrow\infty}\cfrac{|A\cap \{1,...,n\}|}{n}$$ provided this limit exists.  Since classically subsets of $\mathbb{N}$ are in bijective correspondence with sequences of $0$s and $1$s, this notion of density still works for such sequences and can provide an account of frequentist probability. Here we study such a notion of frequentist probability in a constructive framework and we introduce some probabilistic forms of the limited principles of omniscience ($\mathbf{LPO}$) and show that they are either true, false or equivalent to $\mathbf{LPO}$ itself. Our treatment is informal \`a la Bishop (see \cite{BB85}), but it can be formalized in the extensional level of the Minimalist Foundation (see \cite{mtt}, \cite{m09}) with the addition of the axiom of unique choice. 
\section{A constructive account of natural density}
\subsection{Potential events}
A \emph{potential event} $\mathsf{e}$ is a sequence of $0$s and $1$s
$$\mathsf{e}(n)\in \{0,1\}\;[\,n\in \mathbb{N}^{+}]$$

Potential events form a set with extensional equality, that is two potential events $\mathsf{e}$ and $\mathsf{e'}$ are equal (we write $\mathsf{e}=_{\mathcal{P}}\mathsf{e}'$) if  $\mathsf{e}(n)=\mathsf{e}'(n)\;[\,n\in \mathbb{N}^{+}]$.
This set, which we denote with $\mathcal{ P}$, can be endowed with a structure of Boolean algebra as follows:
\begin{enumerate}
\item the bottom $\bot$ is $\lambda n.0$;
\item the top $\top$ is $\lambda n.1$;
\item the conjunction $\mathsf{e}\wedge \mathsf{e}'$ of $\mathsf{e}$ and $ \mathsf{e}'$ in $\mathcal{P}$ is $\lambda n.(\mathsf{e}(n)  \mathsf{e}'(n))$;
\item the disjunction of $\mathsf{e}\vee \mathsf{e}'$ of $\mathsf{e}$ and $ \mathsf{e}'$ in $\mathcal{P}$ is $\lambda n.(\mathsf{e}(n)+  \mathsf{e}'(n)-\mathsf{e}(n)  \mathsf{e}'(n))$;
\item the negation $\neg\mathsf{e}$ of $\mathsf{e}\in \mathcal{P}$ is $\lambda n.(1-\mathsf{e}(n))$;
\item for every $\mathsf{e},\mathsf{e}'\in \mathcal{P}$, $\mathsf{e}\leq \mathsf{e}'$ if and only if $\mathsf{e}(n)\leq \mathsf{e}'(n)\;[\,n\in \mathbb{N}^{+}]$.
\end{enumerate}

Potential events should be understood as sequences of outcomes (``success'' or ``failure'') in a sequence of iterated trials. 

For every potential event $\mathsf{e}$, one can define a sequence of rational numbers
$$\Phi(\mathsf{e})(n)\in \mathbb{Q}\,[\,n\in \mathbb{N}^{+}]$$
taking $\Phi(\mathsf{e})(n)$ to be $\frac{\sum_{i=1}^{n}\mathsf{e}(i)}{n}$. This sequence is called the sequence of \emph{rates of success} (or \emph{frequencies}) of the potential event $\mathsf{e}$.

\subsection{Actual events}
Actual events are those potential events for which the sequence of rates of success can be shown constructively to be convergent.
\begin{definition}
An \emph{actual event} is a pair $(\mathsf{e}, \gamma)$ where $\mathsf{e}$ is a potential event and $\gamma$ is a strictly increasing sequence of natural  numbers such that 
$$\left|\Phi(\mathsf{e})(\gamma(n)+i)-\Phi(\mathsf{e})(\gamma(n)+j)\right|\leq\frac{1}{n}\;\Big[\,n\in \mathbb{N}^{+},\ i\in \mathbb{N},\ j\in \mathbb{N}\Big]$$
We denote with $\widetilde{\mathcal{ A}}$ the set of actual events and two actual events $(\mathsf{e}, \gamma)$ and $( \mathsf{e}', \gamma')$ are \emph{equal}, $(\mathsf{e},\gamma)=_{\widetilde{\mathcal{A}}}(\mathsf{e}',\gamma')$, if $\mathsf{e}=_{\mathcal{P}} \mathsf{e}'$. \end{definition}

Let us now recall the definition of Bishop real numbers from \cite{BB85}.

\begin{definition}
A \emph{Bishop real} $x$ is a sequence $x(n)\in \mathbb{Q}\;[n\in \mathbb{N}^{+}]$ of rational numbers 
such that
$$\left|x(n)-x(m)\right|\leq\cfrac{1}{n}+\cfrac{1}{m}\;\Big[\,n\in \mathbb{N}^{+},\ m\in \mathbb{N}^{+}\Big]$$
Two Bishop reals $x$ and $y$ are \emph{equal} if $\left|x(n)-y(n)\right|\leq\frac{2}{n}\;[\,n\in \mathbb{N}^{+}]$.
The set of Bishop reals is denoted with $\mathbb{R}$ and the equality between Bishop reals with $=_{\mathbb{R}}$.
\end{definition}

We now show that one can define a notion of probability on actual events.
\begin{proposition}
If $(\mathsf{e},\gamma)$ is an actual event, then $\Phi(\mathsf{e})\circ \gamma$ is a Bishop real number. 
\end{proposition}
\begin{proof} Suppose $m\leq n\in \mathbb{N}^{+}$. Then $\gamma(m)\leq\gamma(n)$ and hence
$$\left|\Phi(\mathsf{e}\circ\gamma)(n)-\Phi(\mathsf{e}\circ\gamma)(m)\right|=\left|\Phi(\mathsf{e})(\gamma(n))-\Phi(\mathsf{e})(\gamma(m))\right|\leq\frac{1}{m}<\frac{1}{n}+\frac{1}{m}.$$
This implies that $\Phi(\mathsf{e})\circ \gamma$ is a Bishop real. \end{proof}

\begin{proposition}
If $(\mathsf{e},\gamma)$ and $(\mathsf{e}',\gamma')$ are equal actual events, then $\Phi(\mathsf{e})\circ \gamma$ and $\Phi(\mathsf{e}')\circ \gamma'$ are equal Bishop reals.\end{proposition}

\begin{proof} 
Suppose $(\mathsf{e},\gamma)$ and $(\mathsf{e}',\gamma')$ are equal actual events. Since, for every positive natural number $n$, $\gamma(n)\leq \gamma'(n)$ or $\gamma'(n)\leq \gamma(n)$ and $\mathsf{e}=_{\mathcal{P}}\mathsf{e}'$, then 
$$\left|(\Phi(\mathsf{e})\circ \gamma)(n)-(\Phi(\mathsf{e}')\circ \gamma')(n)\right|=\left|\Phi(\mathsf{e})(\gamma(n))-\Phi(\mathsf{e})(\gamma'(n))\right|\leq \frac{1}{n}<\frac{2}{n}\;\Big[\,n\in \mathbb{N}^{+}\,\Big].$$
From this it follows that $\Phi(\mathsf{e})\circ \gamma$ and $\Phi(\mathsf{e}')\circ \gamma'$ are equal Bishop reals. \end{proof}

As a consequence of the previous two propositions we can give the following 
\begin{definition}
The function $\mathbb{P}:\widetilde{\mathcal{ A}}\rightarrow \mathbb{R}$ is defined as follows: if $(\mathsf{e},\gamma)$ is an actual event, then $\mathbb{P}(\mathsf{e},\gamma):=\Phi(\mathsf{e})\circ \gamma$.
\end{definition}

\subsection{Properties of actual events}
We show here some basic properties of $\mathbb{P}$. However, before proceeding, we need to notice that the following trivial result holds:
\begin{lemma}\label{greater}
If $(\mathsf{e},\gamma)$ is an actual event and $\gamma'$ is an increasing sequence of positive natural numbers such that
$\gamma(n)\leq \gamma'(n)$ for every $n\in \mathbb{N}^{+}$, 
then $(\mathsf{e},\gamma')$ is an actual event.
\end{lemma}

Obviously impossible events are actual null events, i.e.\ we have \emph{strictness}: 

\begin{proposition}\label{zero}
$(\bot,\lambda n.n)$ is an actual event and $\mathbb{P}(\bot,\lambda n.n)=_{\mathbb{R}}0$.
\end{proposition}
Then we show that actual events and $\mathbb{P}$ satisfy \emph{involution}.

\begin{proposition}If $(\mathsf{e},\gamma)$ is an actual event, then $(\neg \mathsf{e},\gamma)$ is an actual event and $$\mathbb{P}(\neg \mathsf{e},\gamma)=_{\mathbb{R}}1-\mathbb{P}(\mathsf{e},\gamma)\footnote{We follow here the definitions of the operations and relations between reals given in \cite{BB85}.}$$
\end{proposition}
\begin{proof} 
Let $(\mathsf{e},\gamma)$ be an actual event.
It is immediate to see that 

$\Phi(\neg \mathsf{e})(n)=1-\Phi(\mathsf{e})(n)$ for every $n\in \mathbb{N}^{+}$. From this it follows that 
$$\left|\Phi(\neg \mathsf{e})( \gamma(n)+m)-\Phi(\neg \mathsf{e})( \gamma(n)+m')\right|=$$
$$=\left|\Phi(\mathsf{e})( \gamma(n)+m')-\Phi(\mathsf{e})( \gamma(n)+m)\right|\leq \frac{1}{n}\; \Big[\,n\in \mathbb{N}^{+},\ m\in \mathbb{N},\ m'\in \mathbb{N}\,\Big].$$
Hence $(\neg \mathsf{e},\gamma)$ is an actual event.  
Moreover for every $n\in \mathbb{N}^{+}$
$$\left|\mathbb{P}(\neg \mathsf{e},\gamma)(n)-(1-\mathbb{P}(\mathsf{e},\gamma))(n)\right|=\left|1-\mathbb{P}(\mathsf{e},\gamma)(n)-(1-\mathbb{P}(\mathsf{e},\gamma)(2n))\right|=$$
$$=\left|\mathbb{P}(\mathsf{e},\gamma)(2n)-\mathbb{P}(\mathsf{e},\gamma)(n)\right|\leq\frac{1}{n}+\frac{1}{2n}<\frac{2}{n}$$
so $\mathbb{P}(\neg \mathsf{e},\gamma)=_{\mathbb{R}}1-\mathbb{P}(\mathsf{e},\gamma)$. \end{proof}

Moreover actual events are closed under \emph{disjunction of incompatible actual events}:

\begin{proposition}\label{disj}
If $(\mathsf{e},\gamma)$ and $(\mathsf{e}',\gamma')$ are actual events with $\mathsf{e}\wedge\mathsf{e}'=_{\mathcal{P}}\bot$ and $\eta:=\lambda n.(\gamma(2n)+\gamma'(2n))$, then the pair $(\mathsf{e}\vee \mathsf{e}', \eta)$ is an actual event.
\end{proposition}
\begin{proof} 
Since $\mathsf{e}\wedge \mathsf{e}'=\bot$, we have that $\mathsf{e}\vee \mathsf{e}'=\mathsf{e}+\mathsf{e}'$. 
Moreover $\eta$ is clearly stricly increasing, as $\gamma$ and $\gamma'$ are so.
Moreover, for every $n\in \mathbb{N}^{+}$, it is immediate to see that $\Phi(\mathsf{e}+\mathsf{e}',n)=\Phi(\mathsf{e},n)+\Phi(\mathsf{e}',n)$.

In particular for every $n\in \mathbb{N}^{+}$, $i,j\in \mathbb{N}$
$$\left|\Phi(\mathsf{e}+\mathsf{e}')(\eta(n)+i)-\Phi(\mathsf{e}+\mathsf{e}')(\eta(n)+j)\right|\leq$$
$$\leq\left|\Phi(\mathsf{e})(\gamma(2n)+\gamma'(2n)+i)-\Phi(\mathsf{e})(\gamma(2n)+\gamma'(2n)+j)\right|+$$
$$\left|\Phi(\mathsf{e}')(\gamma(2n)+\gamma'(2n)+i)-\Phi(\mathsf{e}')(\gamma(2n)+\gamma'(2n)+j)\right|\leq\frac{1}{2n}+\frac{1}{2n}=\frac{1}{n}.$$

Hence $(\mathsf{e}\vee \mathsf{e}',\eta)$ is an actual event.  \end{proof}

We now show that the set of actual events with \emph{null} probability is \emph{downward closed}:
\begin{proposition}
Let $(\mathsf{e},\gamma)$ be an actual event with $\mathbb{P}(\mathsf{e},\gamma)=_{\mathbb{R}}0$ and let $\mathsf{e}'$ be a potential event such that $\mathsf{e}'\leq \mathsf{e}$, then $(\mathsf{e}',\lambda n.\gamma(6n))$ is an actual event.
\end{proposition}
\begin{proof} For every $n\in \mathbb{N}^{+}$ and $i,j\in \mathbb{N}$
$$\left|\Phi(\mathsf{e}')(\gamma(6n)+i)-\Phi(\mathsf{e}')(\gamma(6n)+j)\right|\leq \left|\Phi(\mathsf{e}')(\gamma(6n)+i)\right|+\left|\Phi(\mathsf{e}')(\gamma(6n)+j)\right|\leq $$
$$\leq \left|\Phi(\mathsf{e})(\gamma(6n)+i)\right|+\left|\Phi(\mathsf{e})(\gamma(6n)+j)\right|\leq $$
$$ \left|\Phi(\mathsf{e})(\gamma(6n)+i)-\Phi(\mathsf{e})(\gamma(6n))\right|+\left|\Phi(\mathsf{e})(\gamma(6n)+j)-\Phi(\mathsf{e})(\gamma(6n))\right|+2\left|\Phi(\mathsf{e})(\gamma(6n))\right|$$
$$\leq \frac{1}{6n}+\frac{1}{6n}+2\frac{2}{6n}=\frac{1}{n}$$
by using the fact that $(\mathsf{e},\gamma)$ is an actual event and $\mathbb{P}((\mathsf{e},\gamma))=_{\mathbb{R}}0$.
Hence $(\mathsf{e}',\lambda n.\gamma(6n))$ is an actual event.
\end{proof}
Moreover $\mathbb{P}$ satisfies \emph{monotonicity}:
\begin{proposition} If $(\mathsf{e},\gamma)$ and $(\mathsf{e}',\gamma')$ are actual events and $\mathsf{e}\leq \mathsf{e}'$, then 
$$\mathbb{P}(\mathsf{e},\gamma)\leq \mathbb{P}(\mathsf{e}',\gamma').$$
\end{proposition}
\begin{proof} 
First of all if we take $\eta$ to be the sequence $\lambda n.(\gamma(n)+\gamma'(n))$, then $(\mathsf{e},\eta)$ and $(\mathsf{e}',\eta)$ are actual events  equal to $(\mathsf{e},\gamma)$ and $(\mathsf{e}',\gamma')$ respectively as a consequence of proposition \ref{greater}. In order to show that $\mathbb{P}(\mathsf{e},\eta)\leq \mathbb{P}(\mathsf{e}',\eta)$, we must prove that
$(\mathbb{P}(\mathsf{e},\eta)-\mathbb{P}(\mathsf{e}',\eta))(n)\leq \frac{1}{n}$ for every $n\in \mathbb{N}^{+}$.
But $(\mathbb{P}(\mathsf{e},\eta)-\mathbb{P}(\mathsf{e}',\eta))(n)$ is equal to 
$$(\mathbb{P}(\mathsf{e},\eta))(2n)-(\mathbb{P}(\mathsf{e}',\eta))(2n)=\Phi(\mathsf{e})( \eta(2n))-\Phi(\mathsf{e}')(\eta(2n))$$
Hence we must prove that $\Phi(\mathsf{e})(\eta(2n))\leq \Phi(\mathsf{e}')(\eta(2n))+\frac{1}{n}\; [\,n\in \mathbb{N}^{+}\,]$.
But this is trivially true, because from $\mathsf{e}\leq \mathsf{e}'$ we can deduce that $\Phi(\mathsf{e})(n)\leq \Phi(\mathsf{e}')(n)$ for every positive natural number $n$.
\end{proof}
Finally we prove that $\mathbb{P}$ satisfies a form of \emph{modularity}:
\begin{proposition}\label{modulo}
If $(\mathsf{e},\alpha)$, $(\mathsf{e}',\beta)$, $(\mathsf{e}\wedge \mathsf{e}', \gamma)$, $(\mathsf{e}\vee \mathsf{e}', \delta)$ are actual events, then 
$$\mathbb{P}(\mathsf{e}\vee \mathsf{e}', \delta)+\mathbb{P}(\mathsf{e}\wedge \mathsf{e}', \gamma)=_{\mathbb{R}}\mathbb{P}(\mathsf{e},\alpha)+\mathbb{P}(\mathsf{e}',\beta).$$
\end{proposition}
\begin{proof} 
Using proposition \ref{greater}, if $\varepsilon:=\lambda n.(\alpha(n)+\beta(n)+\gamma(n)+\delta(n))$, then
\begin{enumerate} 
\item $\mathbb{P}(\mathsf{e}\vee \mathsf{e}', \delta)+\mathbb{P}(\mathsf{e}\wedge \mathsf{e}', \gamma)=_{\mathbb{R}}\mathbb{P}(\mathsf{e}\vee \mathsf{e}', \varepsilon)+\mathbb{P}(\mathsf{e}\wedge \mathsf{e}', \varepsilon)$
\item $\mathbb{P}(\mathsf{e},\alpha)+\mathbb{P}(\mathsf{e}',\beta)=_{\mathbb{R}}\mathbb{P}(\mathsf{e},\varepsilon)+\mathbb{P}(\mathsf{e}',\varepsilon)$
\end{enumerate}
So we must prove that  $\mathbb{P}(\mathsf{e},\varepsilon)+\mathbb{P}(\mathsf{e}',\varepsilon)=\mathbb{P}(\mathsf{e}\vee \mathsf{e}', \varepsilon)+\mathbb{P}(\mathsf{e}\wedge \mathsf{e}', \varepsilon)$. However this is immediate since $\mathsf{e}(n)+\mathsf{e}'(n)=(\mathsf{e}\vee \mathsf{e}')(n)+(\mathsf{e}\wedge \mathsf{e}')(n)$ for every $n\in \mathbb{N}^{+}$ 
from which it immediately follows that $\Phi(\mathsf{e},n)+\Phi(\mathsf{e}',n)=\Phi(\mathsf{e}\wedge \mathsf{e}',n)+\Phi(\mathsf{e}\vee \mathsf{e}',n)$ for every $n\in \mathbb{N}^{+}$.
\end{proof}
\subsection{Regular events}
Among potential events, there are some events which can be considered sort of ``deterministic'', since their sequence of outcomes of trials have a periodic behaviour, up to a possible finite number of accidental errors in the recording of the results. These events are here called \emph{regular}:
\begin{definition}
Let $\underline{\alpha}$ and $\underline{\pi}$ be two finite lists of elements of $\{0,1\}$ with length $\ell(\underline{\alpha})$ and $\ell(\underline{\pi})>0$, respectively. 
We define the potential event $\left\|\underline{\alpha},\underline{\pi}\right\|$ as follows
\begin{equation}\notag
\begin{cases}
\left\|\underline{\alpha},\underline{\pi}\right\|(i):=\alpha_{i}\textrm{ if }i\in \mathbb{N}^{+}, i\leq \ell(\underline{\alpha})\\
\left\|\underline{\alpha},\underline{\pi}\right\|(i):=\pi_{\mathsf{rm}(i-\ell(\underline{\alpha})-1,\ell(\underline{\pi}))+1}\;\textrm{ if }i>\ell(\underline{\alpha})\\
\end{cases}
\end{equation}
where $\alpha_{i}$ and $\pi_{i}$ denote the $i$th component of $\underline{\alpha}$ and $\underline{\pi}$, respectively, and where for all $a\in \mathbb{N}$ and $b\in \mathbb{N}^{+}$, $\mathsf{rm}(a,b)$ is the remainder of $a$ divided by $b$. 
The potential events of this form are called \emph{regular}.\end{definition}

The goal of the following propositions is to show that regular events are actual. 

First we show that regular events without errors (that is, events for which $\underline{\alpha}$ is the empty list) are actual events and their probability is equal to the frequency of their period.

\begin{proposition}\label{pi}For every finite non-empty list $\underline{\pi}=[\pi_{1},...,\pi_{m}]$, 
$\left(\left\|[\,],[\underline{\pi}]\right\|,\lambda n. 4nm\right)$ is an actual event and 
$$\mathbb{P}(\left\|[\,],\underline{\pi}\right\|,\lambda n.4nm)=_{\mathbb{R}}\frac{\sum_{k=1}^{m}\pi_{k}}{m}.$$
\end{proposition}
\begin{proof} For every $n\in \mathbb{N}^{+}$ and every $i,j\in \mathbb{N}$ 
{\small
$$\left|\Phi(\left\|[\,],\underline{\pi}\right\|)(4nm+i)-\Phi(\left\|[\,],\underline{\pi}\right\|)(4nm+j)\right|=$$

$$=\left|\frac{\sum_{k=1}^{m}\pi_{k}}{m}\left(\cfrac{(4n+\mathsf{qt}(i,m))m}{4nm+i}\right)+\cfrac{\sum_{k=1}^{\mathsf{rm}(i,m)}\pi_{k}}{4nm+i}-\frac{\sum_{k=1}^{m}\pi_{k}}{m}\left(\cfrac{(4n+\mathsf{qt}(j,m))m}{4nm+j}\right)-\cfrac{\sum_{k=1}^{\mathsf{rm}(j,m)}\pi_{k}}{4nm+j}\right|\leq$$

$$\leq \frac{\sum_{k=1}^{m}\pi_{k}}{m}\left|\cfrac{(4n+\mathsf{qt}(i,m))m}{4nm+i}-\cfrac{(4n+\mathsf{qt}(j,m))m}{4nm+j}\right|+\left|\cfrac{\sum_{k=1}^{\mathsf{rm}(i,m)}\pi_{k}}{4nm+i}+\cfrac{\sum_{k=1}^{\mathsf{rm}(j,m)}\pi_{k}}{4nm+j}\right|\leq$$

$$\leq \left|\cfrac{(4n+\mathsf{qt}(i,m))m}{4nm+i}-\cfrac{(4n+\mathsf{qt}(j,m))m}{4nm+j}\right|+\left|\cfrac{2\sum_{k=1}^{m}\pi_{k}}{4nm}\right|
\leq 1-\cfrac{4n}{4n+1}+\frac{1}{2n}= \frac{1}{4n+1}+\frac{1}{2n}<\frac{1}{n}$$
}
Hence $\left(\left\|[\,],\underline{\pi}\right\|,\lambda n.4nm\right)$ is an actual event.
Moreover for every $i\in \mathbb{N}^{+}$
$$\mathbb{P}(\left\|[\,],\underline{\pi}\right\|,\lambda n.4nm)(i)=\frac{\sum_{k=1}^{m}\pi_{k}}{m}$$ Hence $\mathbb{P}(\left\|[\,],\underline{\pi}\right\|,\lambda n.4nm)=_{\mathbb{R}}\cfrac{\sum_{k=1}^{m}\pi_{k}}{m}.$ \end{proof}

The next step consists in proving that regular events definitely equal to $0$ are actual events and their probability is $0$.

\begin{proposition}\label{alfa} If $\underline{\alpha}$ is a finite list of $0$s and $1$s with length $m$, then 
$\left(\left\|\underline{\alpha},[0]\right\|,\lambda n.2nm\right)$ is an actual event and $\mathbb{P}(\left\|\underline{\alpha},[0]\right\|,\lambda n. 2nm)=_{\mathbb{R}}0$.
\end{proposition}
\begin{proof} 
Suppose $n\in \mathbb{N}^{+}$ and $i,j\in \mathbb{N}$ with $i\geq j$. Then 
$$\left|\Phi\left(\left\|\underline{\alpha},[0]\right\|)(2nm+i\right)-\Phi\left(\left\|\underline{\alpha},[0]\right\|)(2nm+j\right)\right|\leq$$
$$\leq \left|\Phi\left(\left\|\underline{\alpha},[0]\right\|)(2nm+i\right)\right|+\left|\Phi\left(\left\|\underline{\alpha},[0]\right\|)(2nm+j\right)\right|=$$
$$=\cfrac{\alpha_{1}+...+\alpha_{m}}{2nm+i}+\cfrac{\alpha_{1}+...+\alpha_{m}}{2nm+j}\leq \frac{2m}{2nm}=\frac{1}{n}$$
Hence $\left(\left\|\underline{\alpha},[0]\right\|,\lambda n.2nm\right)$ is an actual event. Moreover for every $n\in \mathbb{N}^{+}$ 
$$(\mathbb{P}(\left\|\underline{\alpha},[0]\right\|,\lambda n.2nm)(n)=\cfrac{\alpha_{1}+...+\alpha_{m}}{2nm}\leq \frac{1}{2n}< \frac{2}{n}$$

Hence $\mathbb{P}(\left\|\underline{\alpha},[0]\right\|,\lambda n.2nm)=_{\mathbb{R}}0$.
\end{proof}

Finally we prove that actual events are closed under shift to the right of terms and that probability is preserved by these shifts.

\begin{definition} If $\mathsf{e}$ is a potential event, then $\mathsf{e}^{+}$ is the potential event defined by
\begin{equation}\notag
\begin{cases}
\mathsf{e}^{+}(1):=0\\
\mathsf{e}^{+}(n+1):=\mathsf{e}(n)\textrm{ for }n\in \mathbb{N}^{+}
\end{cases}
\end{equation}
\end{definition}

\begin{proposition}\label{piu}
If $(\mathsf{e},\gamma)$ is an actual event, then $(\mathsf{e}^{+},\lambda n. (\gamma(3n)+1))$ is an actual event and 
$\mathbb{P}(\mathsf{e}^{+},\lambda n.(\gamma(3n)+1))=_{\mathbb{R}}\mathbb{P}(\mathsf{e},\gamma)$.
\end{proposition}
\begin{proof} If $n\in \mathbb{N}^{+}$ and $i,j\in \mathbb{N}$, then $\left|\Phi(\mathsf{e}^{+})(\gamma(3n)+1+i)-\Phi(\mathsf{e}^{+})(\gamma(3n)+1+j)\right|$ is less or equal than the sum of
\begin{enumerate}
\item $ \left|\Phi(\mathsf{e}^{+})(\gamma(3n)+1+i)-\Phi(\mathsf{e})(\gamma(3n)+i)\right|$, 
\item $\left|\Phi(\mathsf{e})(\gamma(3n)+i)-\Phi(\mathsf{e})(\gamma(3n)+j)\right|$ and 
\item $\left|\Phi(\mathsf{e})(\gamma(3n)+j)-\Phi(\mathsf{e}^{+})(\gamma(3n)+1+j)\right|$.
\end{enumerate}
However this sum is equal to 
$$\frac{\Phi(\mathsf{e})(\gamma(3n)+i)}{\gamma(3n)+i+1}+\left|\Phi(\mathsf{e})(\gamma(3n)+i)-\Phi(\mathsf{e})(\gamma(3n)+j)\right|+ \frac{\Phi(\mathsf{e})(\gamma(3n)+j)}{\gamma(3n)+j+1}$$
$$\leq \frac{1}{\gamma(3n)+i+1}+\frac{1}{3n}+\frac{1}{\gamma(3n)+j+1}\leq 3\frac{1}{3n}=\frac{1}{n} $$
Hence $(\mathsf{e}^{+},\lambda n. (\gamma(3n)+1))$ is an actual event. Moreover 
$$\left|\Phi(\mathsf{e}^{+})(\gamma(3n)+1)-\Phi(\mathsf{e})(\gamma(n))\right|=\left|\Phi(\mathsf{e})(\gamma(3n))\frac{\gamma(3n)}{\gamma(3n)+1}-\Phi(\mathsf{e})(\gamma(n))\right|\leq $$
$$\left|\Phi(\mathsf{e})(\gamma(3n))-\Phi(\mathsf{e})(\gamma(n))\right|+\left|\Phi(\mathsf{e})(\gamma(3n))\frac{1}{\gamma(3n)+1}\right|\leq \frac{1}{n}+\frac{1}{3n+1}<\frac{2}{n}$$
Hence $\mathbb{P}(\mathsf{e}^{+},\lambda n.(\gamma(3n)+1))=_{\mathbb{R}}\mathbb{P}(\mathsf{e},\gamma)$. \end{proof}

Putting these results together we obtain the following:
\begin{theorem}
For every regular event $\left\|\underline{\alpha},\underline{\pi}\right\|$, there exists a strictly increasing sequence of natural numbers $\gamma$ such that $(\left\|\underline{\alpha},\underline{\pi}\right\|,\gamma)$ is an actual event and $\mathbb{P}(\left\|\underline{\alpha},\underline{\pi}\right\|,\gamma)=\frac{\sum_{k=1}^{\ell(\underline{\pi})}\pi_{k}}{\ell(\underline{\pi})}$. 
\end{theorem}
\begin{proof}
First of all, notice that $\left\|\underline{\alpha},\underline{\pi}\right\|=_{\mathcal{P}}\left\|\underline{\alpha},[0]\right\|\vee (...(\left\|[\;],\underline{\pi}\right\|)^{+}...)^{+}$ where $^{+}$ is applied $\ell(\underline{\alpha})$ times. 
By proposition \ref{alfa}, we know that there exists $\gamma_{1}$ such that $(\left\|\underline{\alpha},[0]\right\|,\gamma_{1})$ is an actual event and, as a consequence of propositions \ref{pi} and \ref{piu}, there exists $\gamma_{2}$ such that $((...(\left\|[\;],\underline{\pi}\right\|)^{+}...)^{+},\gamma_{2})$ is an actual event. 

Moreover it is clear that $\left\|\underline{\alpha},[0]\right\|\wedge (...(\left\|[\;],\underline{\pi}\right\|)^{+}...)^{+}=_{\mathcal{P}}\bot$, thus, using proposition \ref{disj}, we obtain that there exists  
$\gamma$ such that $(\left\|\underline{\alpha},\underline{\pi}\right\|,\gamma)$ is an actual event.  

Finally, $\mathbb{P}(\left\|\underline{\alpha},\underline{\pi}\right\|,\gamma)=\mathbb{P}(\left\|\underline{\alpha},[0]\right\|,\gamma_{1})+\mathbb{P}((...(\left\|[\;],\underline{\pi}\right\|)^{+}...)^{+},\gamma_{2})-\mathbb{P}(\bot,\lambda n.n)$ by propositions \ref{modulo} and \ref{zero} and, by propositions $\ref{alfa}$, $\ref{piu}$, $\ref{pi}$ and $\ref{zero}$, this is equal to 
$\frac{\sum_{k=1}^{\ell(\underline{\pi})}\pi_{k}}{\ell(\underline{\pi})}$.
\end{proof}

We conclude this section with the following result:

\begin{theorem} 
Regular events form a boolean algebra with the operations inherited by the algebra of potential events.
\end{theorem}
\begin{proof}
First of all $\bot=\left\|[],[0]\right\|$ and $\top=\left\|[],[1]\right\|$. Suppose now that $\left\|\underline{\alpha},\underline{\pi}\right\|$ and $\left\|\underline{\beta},\underline{\psi}\right\|$ are regular events. Then $\neg\left\|\underline{\alpha},\underline{\pi}\right\|:=\left\|\neg\underline{\alpha},\neg\underline{\pi}\right\|$ where $\neg\underline{\alpha}$ and $\neg\underline{\pi}$ are obtained by changing each term $x$ of the finite lists in $1-x$. 
Without loss of generality, we can suppose that $\ell(\underline{\beta})\geq\ell(\underline{\alpha})$. It is immediate to check that 

$$\left\|\underline{\alpha},\underline{\pi}\right\|\wedge\left\|\underline{\beta},\underline{\psi}\right\|=\left\|\underline{\gamma},\underline{\rho}\right\|$$
$$\left\|\underline{\alpha},\underline{\pi}\right\|\vee\left\|\underline{\beta},\underline{\psi}\right\|=\left\|\underline{\gamma'},\underline{\rho'}\right\|$$
where $\ell(\underline{\gamma})=\ell(\underline{\gamma'})=\ell(\underline{\beta})$, $\ell(\underline{\rho})=\ell(\underline{\rho'})=\ell(\underline{\pi})\ell(\underline{\psi})$,  
\begin{equation}\notag
\begin{cases}
\gamma_{i}:=\alpha_{i}\wedge\beta_{i}\textrm{ if }i\leq \ell(\underline{\alpha})\\
\gamma_{i}:=\pi_{\mathsf{rm}(i-\ell(\underline{\alpha})-1,\ell(\underline{\pi}))+1}\wedge \beta_{i}\textrm{ if }\ell(\underline{\alpha})<i\leq \ell(\underline{\beta})
\end{cases}
\end{equation}
\begin{equation}\notag
\begin{cases}
\gamma'_{i}:=\alpha_{i}\vee\beta_{i}\textrm{ if }i\leq \ell(\underline{\alpha})\\
\gamma'_{i}:=\pi_{\mathsf{rm}(i-\ell(\underline{\alpha}-1,\ell(\underline{\pi}))+1}\vee \beta_{i}\textrm{ if }\ell(\underline{\alpha})<i\leq \ell(\underline{\beta})
\end{cases}
\end{equation}
\begin{equation}\notag
\begin{cases}
\rho_{i}:=\pi_{\mathsf{rm}(\ell(\underline{\beta})-\ell(\underline{\alpha})+i-1,\ell(\underline{\pi}))+1}\wedge \psi_{\mathsf{rm}(i-1,\ell(\underline{\psi}))+1}\\
\rho'_{i}:=\pi_{\mathsf{rm}(\ell(\underline{\beta})-\ell(\underline{\alpha})+i-1,\ell(\underline{\pi}))+1}\vee \psi_{\mathsf{rm}(i-1,\ell(\underline{\psi}))+1}\\
 \textnormal{ for every }1\leq i\leq \ell(\underline{\pi})\ell(\underline{\psi}).
 \end{cases}\end{equation}\end{proof}
 
\section{Probabilistic versions of the limited principle of omniscience?}
The limited principle of omniscience, for short $\mathbf{LPO}$, is a non-constructive principle which is weaker than the law of excluded middle (see e.g.\ \cite{bridges1987varieties}) which plays an important role in constructive reverse mathematics (see \cite{Ishihara2006}). It can be formulated in our framework as follows
$$\mathbf{LPO}:(\forall \mathsf{e}\in\mathcal{P})((\forall n\in \mathbb{N}^{+})(\mathsf{e}(n)=0)\vee (\exists n\in \mathbb{N}^{+})(\mathsf{e}(n)=1))$$
Having introduced a notion of probability $\mathbb{P}$ on (some) binary sequences, we can consider some variants of this principle by restricting it to some subset of potential events or by modifying the disjunction in it.

First, we need some abbreviations: 
\begin{enumerate}
\item $\mathsf{List}^{+}(A)$ denotes the set of lists of elements from a set $A$ with positive length;
\item $\mathsf{Incr}(\mathbb{N}^{+},\mathbb{N}^{+})$ denotes the set of increasing sequences of positive natural numbers.
\item $\mathsf{e}\,\varepsilon\,\mathcal{ A}$ means $(\exists \gamma\in \mathsf{Incr}(\mathbb{N}^{+},\mathbb{N}^{+}))((\mathsf{e},\gamma)\in\widetilde{\mathcal{ A}})$;
\item $\varphi(\mathbb{P}[\mathsf{e}])$ means $\mathsf{e}\,\varepsilon\,\mathcal{ A}\wedge (\forall \gamma\in \mathsf{Incr}(\mathbb{N}^{+},\mathbb{N}^{+}))((\mathsf{e},\gamma)\in\widetilde{\mathcal{ A}}\rightarrow \varphi(\mathbb{P}(\mathsf{e},\gamma)))$ whenever $\varphi(x)$ is a proposition depending on a real number $x$;
\item $\mathsf{e}\,\varepsilon\,\mathcal{ R}$ means $(\exists {\underline \alpha}\in \mathsf{List}(\{0,1\}))(\exists {\underline \pi}\in \mathsf{List}^{+}(\{0,1\}))(\mathsf{e}=_{\mathcal{P}}\left\|\underline{\alpha},\underline{\pi}\right\|)$;
\item $\mathsf{e}\,\varepsilon\,\mathcal{ N}$ means 
$\mathbb{P}[\mathsf{e}]=_{\mathbb{R}}0$.
\end{enumerate}
In such a way we can define subsets $\mathcal{A}$, $\mathcal{R}$ and $\mathcal{N}$ of $\mathcal{P}$ consisting of actual, regular and null events, respectively.

Before proceeding, let us prove a simple, but very important fact.
\begin{lemma}\label{proproposition}
Let $\mathsf{e}\in \mathcal{ P}$;
\begin{enumerate}
\item if $\mathbb{P}[\mathsf{e}]>0$, then $(\exists n\in \mathbb{N}^{+})(\mathsf{e}(n)=1)$;
\item if $(\forall n\in \mathbb{N}^{+})(\mathsf{e}(n)=0)$, then $\mathbb{P}[\mathsf{e}]=_{\mathbb{R}}0$.
\end{enumerate}
\end{lemma}
\begin{proof}\emph{2.} is immediate since $(\forall n\in \mathbb{N}^{+})(\mathsf{e}(n)=0)$ means exactly $\mathsf{e}=_{\mathcal{P}}\bot$.
 
Suppose now that $(\mathsf{e},\gamma)\in \widetilde{\mathcal{ A}}$ and $\mathbb{P}(\mathsf{e},\gamma)>0$. By definition of $>$ between real numbers, there exists $m\in \mathbb{N}$ such that $\mathbb{P}(\mathsf{e},\gamma)(m)>\frac{1}{m}$, i.e.\, 
$$\sum_{i=1}^{\gamma(m)}m\mathsf{e}(i)>\gamma(m)>0$$
Hence there exists $n$ such that $1\leq n\leq \gamma(m)$ and $\mathsf{e}(n)=1$.  

\noindent Thus $(\exists n\in \mathbb{N}^{+})(\mathsf{e}(n)=1)$ and \emph{1.} is proved. \end{proof}

We can now introduce a family of ``probabilistic'' versions of $\mathbf{LPO}$.

\begin{definition}For every subset $\mathcal{ E}$ of $\mathcal{ P}$ we define the following principles: 
\begin{enumerate}
\item {\bf LPO}$[\mathcal{ E}]$ $(\forall \mathsf{e}\in \mathcal{ P})\left[\mathsf{e}\,\varepsilon\,\mathcal{ E}\rightarrow(\forall n\in \mathbb{N}^{+})(\mathsf{e}(n)=0)\vee(\exists n\in \mathbb{N}^{+})(\mathsf{e}(n)=1)\right]$
\item $\mathbb{P}$-{\bf LPO}$[\mathcal{ E}]$ $(\forall \mathsf{e}\in \mathcal{ P})\left[\mathsf{e}\,\varepsilon\,\mathcal{ E}\rightarrow\mathbb{P}[\mathsf{e}]=_{\mathbb{R}}0\vee(\exists n\in \mathbb{N}^{+})(\mathsf{e}(n)=1)\right]$
\item $\mathbb{PP}$-{\bf LPO}$[\mathcal{ E}]$ $(\forall \mathsf{e}\in \mathcal{ P})\left[\mathsf{e}\,\varepsilon\,\mathcal{ E}\rightarrow\mathbb{P}[\mathsf{e}]=_{\mathbb{R}}0\vee\mathbb{P}[\mathsf{e}]>0\right]$
\end{enumerate}
\end{definition}
\noindent Notice that $\mathbf{LPO}[\mathcal{ P}]$ is equivalent to $\mathbf{LPO}$, and that $\mathbb{PP}$-$\mathbf{LPO}[\mathcal{P}]$ implies $\mathcal{P}=\mathcal{A}$. 

Let us now study the relation between these probabilistic versions of $\mathbf{LPO}$ when $\mathcal{E}$ is $\mathcal{P}$, $\mathcal{A}$, $\mathcal{N}$ or $\mathcal{R}$. 
As a direct consequence of lemma \ref{proproposition} we have:
\begin{lemma}\label{propPP}If $\mathcal{ E}$ is a subset of $\mathcal{ P}$, then
\begin{enumerate}
\item $\mathbb{P}\mathbb{P}$-$\mathbf{LPO}[\mathcal{ E}]\rightarrow \mathbb{P}$-$\mathbf{LPO}[\mathcal{ E}]$
\item $\mathbf{LPO}[\mathcal{ E}]\rightarrow \mathbb{P}$-$\mathbf{LPO}[\mathcal{ E}]$
\end{enumerate}
\end{lemma}
Moreover, as a trivial consequence of the definitions we have that
\begin{lemma}\label{propPPP}If $\mathcal{ E}\subseteq \mathcal{ E}'$ are subsets of $\mathcal{ P}$, then
\begin{enumerate}
\item $\mathbf{LPO}[\mathcal{ E}']\rightarrow \mathbf{LPO}[\mathcal{ E}]$
\item $\mathbb{P}$-$\mathbf{LPO}[\mathcal{ E}']\rightarrow \mathbb{P}$-$\mathbf{LPO}[\mathcal{ E}]$
\item $\mathbb{PP}$-$\mathbf{LPO}[\mathcal{ E}']\rightarrow \mathbb{PP}$-$\mathbf{LPO}[\mathcal{ E}]$
\end{enumerate}
\end{lemma}

Some of these principles can be proven constructively, while one is constructively false:
\begin{lemma}The following hold:
\begin{enumerate}
\item $\mathbf{LPO}[\mathcal{ R}]$
\item $\mathbb{P}$-$\mathbf{LPO}[\mathcal{ R}]$
\item $\mathbb{PP}$-$\mathbf{LPO}[\mathcal{ R}]$
\item $\mathbb{PP}$-$\mathbf{LPO}[\mathcal{ N}]$
\item $\mathbb{P}$-$\mathbf{LPO}[\mathcal{ N}]$
\item $\neg\mathbb{PP}$-$\mathbf{LPO}[\mathcal{ P}]$
\end{enumerate}
\end{lemma}
\begin{proof}
Let us prove one by one the statements above.
\begin{enumerate}
\item $\mathbf{LPO}[\mathcal{ R}]$ holds as it reduces to control, for each $\mathsf{e}=_{\mathcal{P}}\left\|\underline{\alpha},\underline{\pi}\right\|$ in $\mathcal{ R}$, a finite number of entries (those in $\underline{\alpha}$ and $\underline{\pi}$). 
\item This is a consequence of \emph{1.\ }and lemma \ref{propPP}. 
\item If $\mathsf{e}=_{\mathcal{P}}\left\|\underline{\alpha},\underline{\pi}\right\|$ in $\mathcal{R}$, then $\mathbb{P}[\mathsf{e}]=_{\mathbb{R}}\cfrac{\sum_{j=1}^{\ell(\underline{\pi})}\pi_{j}}{\ell(\underline{\pi})}\,\varepsilon\, \mathbb{Q}$ and we can thus decide whether $\mathbb{P}[\mathsf{e}]=_{\mathbb{R}}0$ or $\mathbb{P}[\mathsf{e}]>0$.
\item This is obvious, since $\mathsf{e}\,\varepsilon\,\mathcal{ N}$ means, by definition, that $\mathbb{P}[\mathsf{e}]=_{\mathbb{R}}0$.
\item This is a consequence of \emph{4.\ }and lemma \ref{propPP}.
\item There exist potential events $\mathsf{e}$ for which there is no $\gamma$ such that $(\mathsf{e},\gamma)\in \widetilde{\mathcal{ A}}$. Consider for example the sequence 
$10110011110000...$ which alternates a group of $2^{n}$ ones and a group of $2^{n}$ zeros increasing $n$ at each step.
\end{enumerate}
\end{proof}
After these first three lemmas, the situation about the remaining variants can be summarize as follows
{\small 
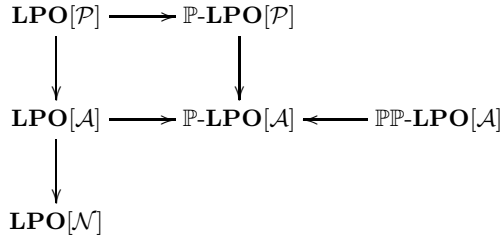
$$\xymatrix{
\mathbf{LPO}[\mathcal{ P}]\ar[r]\ar[d]		&\mathbb{P}\textnormal{-}\mathbf{LPO}[\mathcal{ P}]\ar[d]		\\
\mathbf{LPO}[\mathcal{ A}]\ar[r]\ar[d]			&\mathbb{P}\textnormal{-}\mathbf{LPO}[\mathcal{ A}]	&\mathbb{PP}\textnormal{-}\mathbf{LPO}[\mathcal{ A}]\ar[l]	\\
\mathbf{LPO}[\mathcal{ N}]	&\\
}$$
}

However, the left side of the previous diagram collapses:

\begin{lemma}
$\mathbf{LPO}[\mathcal{N}]$ implies $\mathbf{LPO}[\mathcal{P}]$. Hence $\mathbf{LPO}[\mathcal{P}]$, $\mathbf{LPO}[\mathcal{A}]$ and $\mathbf{LPO}[\mathcal{N}]$ are all equivalent.
\end{lemma} 
\begin{proof} Let $\mathsf{e}$ be a potential event. We can define an increasing sequence $\widetilde{\mathsf{e}}$ in $\mathcal{P}$ as follows:
\begin{equation} \notag
\begin{cases}
\widetilde{\mathsf{e}}(1):=\mathsf{e}(1)\\
\widetilde{\mathsf{e}}(n+1):=\widetilde{\mathsf{e}}(n)\vee\mathsf{e}(n+1)\,[n\in \mathbb{N}^{+}]\\
\end{cases}
\end{equation}
Clearly, $(\exists n\in \mathbb{N}^{+})(\mathsf{e}(n)=1)$ if and only if $(\exists n\in \mathbb{N}^{+})(\widetilde{\mathsf{e}}(n)=1)$, and moreover $(\forall n\in \mathbb{N}^{+})(\mathsf{e}(n)=0)$ if and only if $(\forall n\in \mathbb{N}^{+})(\widetilde{\mathsf{e}}(n)=0)$.

Using $\widetilde{\mathsf{e}}$ we define another potential event $\alpha(\widetilde{\mathsf{e}})$ as follows:
\begin{equation}\notag
\begin{cases}
\alpha(\widetilde{\mathsf{e}})(1):=\widetilde{\mathsf{e}}(1)\\
\alpha(\widetilde{\mathsf{e}})(n+1):=\widetilde{\mathsf{e}}(n+1)-\widetilde{\mathsf{e}}(n)\,[n\in \mathbb{N}^{+}]\\
\end{cases}
\end{equation}
Since $\widetilde{e}$ is increasing, $\alpha(\widetilde{\mathsf{e}})$ takes value $1$ for at most one input. In particular $\alpha(\widetilde{\mathsf{e}})$ is in $\mathcal{N}$ and moreover $(\exists n\in \mathbb{N}^{+})(\widetilde{\mathsf{e}}(n)=1)$ if and only if $(\exists n\in \mathbb{N}^{+})(\alpha(\widetilde{\mathsf{e}})(n)=1)$ and $(\forall n\in \mathbb{N}^{+})(\widetilde{\mathsf{e}}(n)=0)$ if and only if $(\forall n\in \mathbb{N}^{+})(\alpha(\widetilde{\mathsf{e}})(n)=0)$. 

If $\mathbf{LPO}[\mathcal{N}]$ holds, then $(\exists n\in \mathbb{N}^{+})(\alpha(\widetilde{\mathsf{e}})(n)=1)\vee (\forall n\in \mathbb{N}^{+})(\alpha(\widetilde{\mathsf{e}})(n)=0)$, which is equivalent to $(\exists n\in \mathbb{N}^{+})(\mathsf{e}(n)=1)\vee (\forall n\in \mathbb{N}^{+})(\mathsf{e}(n)=0)$. Thus from $\mathbf{LPO}[\mathcal{N}]$ follows $\mathbf{LPO}[\mathcal{P}]$. \end{proof}
After the proof of this lemma the situation is the following:
{\small$$\xymatrix{
\mathbf{LPO}[\mathcal{ P}]\equiv\mathbf{LPO}[\mathcal{ A}]\equiv \mathbf{LPO}[\mathcal{ N}]\ar[d]		&	\mathbb{PP}\textnormal{-}\mathbf{LPO}[\mathcal{ A}]\ar[d]	\\
\mathbb{P}\textnormal{-}\mathbf{LPO}[\mathcal{ P}]\ar[r]		&\mathbb{P}\textnormal{-}\mathbf{LPO}[\mathcal{ A}]				\\
}$$}
\begin{lemma}
$\mathbb{P}$-$\mathbf{LPO}[\mathcal{A}]$ implies $\mathbf{LPO}[\mathcal{P}]$, thus $\mathbb{P}$-$\mathbf{LPO}[\mathcal{A}]$, $\mathbb{P}$-$\mathbf{LPO}[\mathcal{P}]$ and $\mathbf{LPO}[\mathcal{P}]$ are equivalent.
\end{lemma}
\begin{proof}
Let $\mathsf{e}$ be a potential event. It is well known that there exist two primitive recursive functions $\varphi,\psi:\mathbb{N}^{+}\rightarrow \mathbb{N}$ such that, for every $n\in \mathbb{N}^{+}$, $\varphi(n)$ and $\psi(n)$ are the unique natural numbers such that $n=2^{\varphi(n)}(2\psi(n)+1)$.

Let $\widehat{\mathsf{e}}$ be the potential event defined by 
$$\widehat{\mathsf{e}}(n):=\mathsf{e}(\varphi(n)+1)\,[n\in \mathbb{N}^{+}]$$
We have that  $(\exists n\in \mathbb{N}^{+})(\,\widehat{\mathsf{e}}(n)=1\,)$ if and only if $(\exists n\in \mathbb{N}^{+})(\mathsf{e}(n)=1)$.  Moreover $\widehat{\mathsf{e}}\,\varepsilon\,\mathcal{A}$, because $(\widehat{\mathsf{e}},\lambda n.2^{n-1})\in \widetilde{\mathcal{A}}$ and, since for every $n\in \mathbb{N}^{+}$
$$(\mathbb{P}(\widehat{\mathsf{e}},\lambda n.2^{n-1}))(n)=\cfrac{\mathsf{e}(n)}{2^{n-1}}+\sum_{k=1}^{n-1}\frac{\mathsf{e}(k)}{2^{k}}$$
it immediately follows that $\mathbb{P}[\widehat{\mathsf{e}}]=_{\mathbb{R}}\sum_{n=1}^{+\infty}\frac{\mathsf{e}(n)}{2^{n}}$. In particular, $\mathbb{P}[\widehat{\mathsf{e}}]=_{\mathbb{R}}0$ if and only if $(\forall n\in \mathbb{N}^{+})(\mathsf{e}(n)=0)$. 
 
 Assuming $\mathbb{P}$-$\mathbf{LPO}[\mathcal{A}]$, we obtain $\mathbb{P}(\widehat{\mathsf{e}})=0\vee (\exists n\in \mathbb{N}^{+})(\widehat{\mathsf{e}}(n)=1)$ which is equivalent to $(\exists n\in \mathbb{N}^{+})(\mathsf{e}(n)=1)\vee(\forall n\in \mathbb{N}^{+})(\mathsf{e}(n)=0)$.

We have hence shown that $\mathbb{P}$-$\mathbf{LPO}[\mathcal{A}]$ implies $\mathbf{LPO}[\mathcal{P}]$.\end{proof}

After the proof of this lemma the situation is the following:
$$\mathbb{PP}\textnormal{-}\mathbf{LPO}[\mathcal{ A}]\rightarrow \mathbf{LPO}[\mathcal{ P}](\equiv\mathbf{LPO}[\mathcal{ A}]\equiv \mathbf{LPO}[\mathcal{ N}]\equiv \mathbb{P}\textnormal{-}\mathbf{LPO}[\mathcal{ P}]\equiv\mathbb{P}\textnormal{-}\mathbf{LPO}[\mathcal{ A}])$$

However we have the following
\begin{lemma}
$\mathbf{LPO}[\mathcal{P}]$ implies $\mathbb{PP}$-$\mathbf{LPO}[\mathcal{A}]$.
\end{lemma}
\begin{proof}
Assume $\mathbf{LPO}[\mathcal{P}]$. Since we have assumed countable choice, trichotomy for Bishop reals follows (see \cite{bridges1987varieties}) and thus $\mathbb{PP}$-$\mathbf{LPO}[\mathcal{A}]$ holds.
\end{proof}
As consequence of all the lemmas above we have the following
\begin{theorem}
$$\mathbf{LPO}[\mathcal{ P}]\equiv\mathbf{LPO}[\mathcal{ A}]\equiv \mathbf{LPO}[\mathcal{ N}]\equiv \mathbb{P}\textnormal{-}\mathbf{LPO}[\mathcal{ P}]\equiv\mathbb{P}\textnormal{-}\mathbf{LPO}[\mathcal{ A}]\equiv\mathbb{PP}\textnormal{-}\mathbf{LPO}[\mathcal{ A}]$$
\end{theorem}
Thus we have proved that no one of the variations on $\mathbf{LPO}$ which we have introduced is different from a true statement, a false statement or $\mathbf{LPO}$ itself. Thus, from one side the hope of having an authentic ``probabilistic'' variant of $\mathbf{LPO}$ have failed; however the results above provide new statements equivalent to $\mathbf{LPO}$, namely $\mathbf{LPO}[\mathcal{A}]$, $\mathbf{LPO}[\mathcal{N}]$, $\mathbb{P}$-$\mathbf{LPO}[\mathcal{A}]$, $\mathbb{P}$-$\mathbf{LPO}[\mathcal{P}]$ and $\mathbb{PP}$-$\mathbf{LPO}[\mathcal{A}]$.
\subsection*{Acknowledgements}
The author would like to thank Francesco Ciraulo, Hannes Diener, Milly Maietti, Fabio Pasquali and Giovanni Sambin for fruitful discussions about the subject of this paper.
The author thanks EU-MSCA-RISE project 731143 ``Computing with Infinite Data'' (CID) for supporting the research.

\bibliographystyle{plain}
\bibliography{biblioccc}
\end{document}